\documentclass[a4paper,draft]{amsart}
\usepackage[latin1]{inputenc}
\usepackage{amssymb,amsfonts,amsmath,url}
\usepackage[all]{xy}

\usepackage{graphicx}
\usepackage{latexsym,longtable,amscd,booktabs,array}
\usepackage{lscape}
\usepackage{color}
\usepackage[OT2,T1]{fontenc}
\DeclareSymbolFont{cyrletters}{OT2}{wncyr}{m}{n}
\DeclareMathSymbol{\Sha}{\mathalpha}{cyrletters}{"58}

\theoremstyle{plain}
\newtheorem{theorem}{Theorem}

\newtheorem{conjecture}{Conjecture}

\newtheorem{corollary}{Corollary}

\newtheorem{definition}{Definition}
\newtheorem{proposition}{Proposition}

\newcommand{\Q}{\mathbb{Q}}

\newcommand{\Z}{\mathbb{Z}}

\newcommand{{\tors}}{\operatorname{tors}}
\newcommand{\rank}[1]{{\rm rank}\,#1}

\begin{document}

\title[Arithmetic progressions of four squares over quadratic fields]{Arithmetic progressions of four squares over quadratic fields}

\author{Enrique Gonz{\'a}lez-Jim{\'e}nez}
\address{Universidad Aut{\'o}noma de Madrid\\
Departamento de Matem{\'a}ticas \\
28049 Madrid, Spain}
\email{enrique.gonzalez.jimenez@uam.es}

\author{J\"orn Steuding}
\address{W\"urzburg University\\
Institut f\"ur Mathematik\\
97074 W\"urzburg, Germany}
\email{steuding@mathematik.uni-wuerzburg.de}

%\subjclass[2000]{Primary 14G35, 14H45; Secondary 11F11, 11G10}
\date{\today}
\maketitle

\begin{abstract}
Let $d$ be a squarefree integer. Does there exist four squares in arithmetic progression over $\Q(\sqrt{d})$? We shall give a partial answer to this question, depending on the value of $d$. In the affirmative case, we construct explicit arithmetic progressions consisting of four squares over $\Q(\sqrt{d})$.
\end{abstract}

%%%%%%%%%%%%%%%%%%%%%%%%%%%%%%%%%%%%%%%%%%%%%%%%%%
\tableofcontents  % TO BE REMOVED
%%%%%%%%%%%%%%%%%%%%%%%%%%%%%%%%%%%%%%%%%%%%%%%%%%

\section{Introduction}\label{intro}

Non-constant arithmetic progressions consisting of rational squares have been studied since ancient times. While it is not difficult to obtain an arithmetic progression of three rational squares (e.g. $1^2,5^2,7^2$), there are no four distinct rational squares in arithmetic progression as already stated by Fermat and proved by Euler (among others). However, the situation is different over number fields. It is easy to construct four squares in arithmetic progression over a quadratic number field; e.g. $1^2,5^2,7^2,(\sqrt{73})^2$ over $\Q(\sqrt{73})$. It is even possible to find five squares in arithmetic progressions; e.g. $7^2,13^2,17^2,(\sqrt{409})^2,23^2$ over $\Q(\sqrt{409})$. By Faltings' proof of the Mordell conjecture \cite{faltings} in any finite algebraic extension of $\Q$ there can be at most finitely many arithmetic progressions of at least five squares. Recently Xarles \cite{xarles} has proved that six squares in arithmetic progression over quadratic number fields do not exist. The case of length five over quadratic fields has been treated by the first author and Xarles \cite{GJX}.

In this paper we consider the following natural problem:
\vspace{1mm}

\noindent {\it Let $d$ be an squarefree integer. Do there exist four squares in arithmetic progression over $\Q(\sqrt{d})$? In the affirmative case, give an algorithm to construct explicit examples.} \\

For special values of $d$ the following table indicates for which quadratic fields $\Q(\sqrt{d})$ there exist or do not exist non-constant arithmetic progressions of four squares. Here '?' indicates that in this case we do not know whether there exists a non-constant four term arithmetic progression of squares or not.

{\footnotesize
\begin{center}
\begin{table}[!h]

\begin{tabular}{|c||c|c|c|c||c|c|c|c|}
\hline
$p\ge 5$ prime & \multicolumn{8}{|c|}{Is there a non-constant arithmetic progression of four squares over $\Q(\sqrt{d})$?} \\
\hline
$p\bmod\,24$ & $d=p$ & $d=2p$ &$d=3p$ &$d=6p$ & $d=-p$ & $d=-2p$ &$d=-3p$ &$d=-6p$\\
\hline
 $1$ & $?$ & $?$ & $?$ & $?$ & $?$ &  $?$ & $?$ & $?$ \\
\hline
 $5$ & no & $?$ & no & $?$ & $?$ &  $?$ & $?$ & no \\
\hline
 $7$ & no & no &  $?$ & $?$ & no &  $?$ & $?$ & no \\
\hline
 $11$ & $?$ & $?$ & no & $?$ & no &  $?$ & $?$ & $?$ \\
\hline
 $13$ & $?$ & no &  $?$ & $?$ & no &  no & $?$ & no \\
\hline
 $17$ & $?$ & $?$ &  no & $?$ & $?$ &  $?$ & no & no \\
\hline
 $19$ & no & $?$ &  no & $?$ & $?$ & no & $?$ &$?$ \\
\hline
 $23$ & yes & yes &  $?$ & yes & yes &  yes & yes & $?$ \\
\hline
\end{tabular}\\[2mm]
\caption{}\label{table1}
\end{table}
\end{center}
}
\vspace{-8mm}

The proof of the correctness of the table will be given in Section  \ref{section-explicit}. 

We shall show that four squares in arithmetic progression lead to points on an elliptic curve. This approach is not new; it seems to be a folklore result, however, we shall provide our own parametrization of arithmetic progressions of four squares over a number field $k$ by $k$-rational points on the modular curve $X_0(24)$. More precisely, there exists a non-constant arithmetic progression of four squares over a number field $k$ if and only if the Mordell-Weil group of the elliptic curve $X_0(24)$ has positive rank over $k$; in this case, there exist infinitely many arithmetic progressions of four squares. For the specific case of a quadratic number field $\Q(\sqrt{d})$ this characterization reduces our problem to the problem of determining the rank of the quadratic $d$-twists of the underlying elliptic curve.

The above characterization allows us to link our problem with two other problems, both being generalizations of the famous congruent number problem. These problems are related to $\theta$-congruent numbers and Euler's concordant forms.

This paper is organized as follows: Section \ref{section-4squares} is devoted to construct a parametrization of four term arithmetic progressions of squares over a number field $k$ by $k$-rational points of the elliptic curve $X_0(24)$. Here we derive some particular results when $k$ is a quadratic field. In Section \ref{section-congruent} we introduce the notion of $\theta$-congruent numbers and state some well-known results for the case when $\theta$ is either equal to $\pi/3$ or $2\pi/3$, which are the relevant cases with respect to the main problem of this paper. Section \ref{section-concordant} deals with Euler's concordant forms. In Section \ref{section-explicit} we apply results from the previous sections in order to obtain a partial answer to our problem. Here we also give examples of arithmetic progressions of four squares over $\Q(\sqrt{d})$ for all cases $|d|\le 40$ for which such arithmetic progressions do exist. Moreover, we give another construction using pythagorean triples and Thue equations. The last section contains some average results related to our main problem.

\subsection*{Acknowledgements}
Research of the first author was supported in part by grant MTM 2006--10548 (Ministerio de Educaci\'on y Ciencia, Spain) and CCG07--UAM/ESP--1814 (Universidad Aut\'onoma de Madrid -- Comunidad de Madrid, Spain). Research of the second author was supported in part by grant MTM 2006-01859  (Ministerio de Educaci\'on y Ciencia, Spain).

\section{Four squares in arithmetic progression}\label{section-4squares}

Four squares $a^2,b^2,c^2$ and $d^2$ over a field $k$ are in arithmetic progression if and only if $b^2-a^2=c^2-b^2$ and $c^2-b^2=d^2-c^2$. This is equivalent to $[a,b,c,d]\in \mathbb{P}^3(k)$ being in the intersection of the two quadric surfaces $a^2+c^2=2b^2$ and $b^2+d^2=2c^2$, resp. lying on the curve
\begin{equation}\label{eqn1}
C\,:\,\left\{
\begin{array}{l}
a^2+c^2=2b^2,\\
b^2+d^2=2c^2.
\end{array}
\right.
\end{equation}
Therefore, $k$-rational points of $C$ parametrize arithmetic progressions of four squares over $k$. Note that the eight points $[\pm 1,\pm 1,\pm 1,\pm 1 ]$ belong to $C$, however, these points correspond to constant arithmetic progressions. The next step is to compute an explicit equation for $C$. In the generic case the intersection of two quadric surfaces in $\mathbb{P}^3$ gives an elliptic curve and, indeed, this will turn out to be true in our case. For this purpose we are going to compute a Weierstrass equation for $C$. The system of equations (\ref{eqn1}) is equivalent to $a^2+2d^2=3c^2, b^2=2c^2-d^2$. A parametrization (up to sign) of the first conic is given by
$$
(a,d,c)=(2t^2-4t-1,2t^2+2t-1,2t^2+1),\quad t\in k,
$$
where the inverse is given by $t=\frac{d-c}{a-c}$ if $a\neq c$ and $t=-\frac{1}{2}$ if $a=c$. Therefore, if we substitute the values of $a,d$ and $c$ in the second equation, we obtain the quartic equation
\begin{equation}\label{eqn_quartic}
\mathcal{Q}\,:\,b^2=4t^4-8t^3+8t^2+4t+1\,.
\end{equation} 
Our next aim is to find a Weierstrass model for $\mathcal{Q}$. Note that there is only one point at infinity, namely $[0:1:0]$. This point is a node. We denote by $\infty_{1}$ and $\infty_{2}$ the two branches at infinity at the desingularization of $\mathcal{Q}$. A Weierstrass model for $\mathcal{Q}$ is $E:y^2=x(x+3)(x-1)$, where the isomorphism $\phi: \mathcal{Q} \rightarrow E$ is defined by
$$
\phi(P)= \left(\frac{1+b+2t}{2t^2},\frac{1+b+3t+bt+4t^2-2t^3}{2t^3}\right)\quad\mbox{if} \quad P=(t,b)\neq (0,\pm 1),
$$
and $\phi(0,-1)=(-1,2)$, $\phi(0,1)=[0:1:0]$, $\phi(\infty_1)=(1,0)$, $\phi(\infty_2)=(-1,-2)$.  The inverse is defined by
$$
\phi^{-1}(P)= \left(\frac{x+y-1}{x^2-1},\frac{x^3+5x^2+2xy-2y-x+3}{(x^2-1)(x+1)}\right)\quad \mbox{if} \quad P=(x,y),\,x\neq \pm 1.
$$ 
Therefore, by the above construction we have proved

\begin{theorem}\label{parametrization}
Let $k$ be a field of $\mbox{char($k$)}\neq 2,3$, then arithmetic progressions of four squares in $k$ are parametrized by $k$-rational points of the elliptic curve 
$$
E:y^2=x(x+3)(x-1).
$$
This parametrization is as follows:\\
\noindent $\bullet$ Let $[a,b,c,d]\in\mathbb{P}^3(k)$ such that $a^2,b^2,c^2,d^2$ form an arithmetic progression. If $a\neq c$ and $d\neq c$, let $t=\frac{d-c}{a-c}$ and define
$$
P=\left(\frac{1+b+2t}{2t^2},\frac{1+b+3t+bt+4t^2-2t^3}{2t^3}\right),
$$
and 
$$
P=\left\{
\begin{array}{ccl}
[0:1:0] & \mbox{if} &  [a,b,c,d]=[-1,1,1,1],\\ 
(-1,2) & \mbox{if} &  [a,b,c,d]=[-1,-1,1,1],\\
(3,-6)  & \mbox{if} &  [a,b,c,d]=[1,1,1,1],\\
(-3,0) & \mbox{if} &   [a,b,c,d]=[1,-1,1,1],\\
(1,0) & \mbox{if} &  [a,b,c,d]=[-1,1,-1,1],\\
(-1,-2) & \mbox{if} &  [a,b,c,d]=[-1,-1,-1,1].
\end{array}
\right.
$$
Then $P\in E(k)$.\\
\noindent $\bullet$ Let be $P\in E(k)$. If $P=(x,y)$, $x\neq \pm 1$, let $t=\frac{x+y-1}{x^2-1}$ and define
$$
[a,b,c,d]=\left[2t^2-4t-1, \frac{x^3+5x^2+2xy-2y-x+3}{(x^2-1)(x+1)}, 2t^2+1, 2t^2+2t+1\right],
$$
and
$$
[a,b,c,d]=\left\{
\begin{array}{ccl}
{[-1,-1,1,1]} & \mbox{if} & P=(-1,2),\\ 
{[-1,-1,-1,1]} & \mbox{if} & P=(-1,-2),\\ 
{[-1,1,-1,1]} & \mbox{if} & P=(1,0),\\
{[-1,1,1,1]} & \mbox{if} & P=[0:1:0].
\end{array}
\right.
$$
Then $a^2,b^2,c^2,d^2$ form an arithmetic progression in $k$.
\end{theorem} 

It is natural to ask for which number fields $k$ there exist arithmetic progressions of four squares over $k$. In order to investigate this problem we shall make use of the theory of elliptic curves.

\begin{proposition}\label{prop1}
Let $k$ be a number field. Then there exists a non-constant arithmetic progression of four squares over $k$ if and only if $\#E(k)>8$. Furthermore, there exist infinitely many such progressions if and only if $\rank{E(k)}\ne 0$.
\end{proposition}

\begin{proof}
Firstly, note that the points $[\pm 1,\pm 1,\pm 1,\pm 1]$ in $\mathbb{P}^3(k)$ give constant arithmetic progressions. Therefore, using the parametrization $\phi$, it follows that the points $\phi([\pm 1, \pm 1, \pm1, \pm 1])$ belong to $E(k)$, and  this set has cardinality $8$. This concludes the proof.
\end{proof}

As a corollary we obtain

\begin{corollary}
There is no non-constant arithmetic progression of four rational squares.
\end{corollary}

\noindent This statement is due to Fermat, however, the first proof is attributed to Euler who applied Fermat's method of infinite descent. For the sake of completeness and since some of the data that appear will be useful later, we give the short

\vspace{.3cm}

\begin{proof}  Using \verb+SAGE+ \cite{sage} or \verb+MAGMA+ \cite{magma}, one can check that $E$ is the curve \verb+24A1+ in Cremona's tables \cite{cremona}, resp. \verb+24B+ in the Antwerp tables \cite{antwerp}. In other words, $E$ is the modular curve $X_0(24)$. Checking these tables or using one of the above mentioned computer algebra systems, one can prove $E(\Q)\simeq \Z/2\Z\oplus \Z/4\Z$. There are no $\Q$-rational points on $E$ apart those eight points $[\pm 1,\pm 1,\pm 1,\pm 1]$ which correspond to constant arithmetic progressions.
\end{proof}
%\vspace{.3cm}

Next we shall consider the case of quadratic number fields. Here the question translates to study the Mordell-Weil group $E(\Q(\sqrt{d}))$. However, instead to treat the elliptic curve $E$ over $\Q(\sqrt{d})$ directly, we are going to study the quadratic $d$-twist of the elliptic curve $E$ over $\Q$, i.e., the elliptic curve
$$
E^d: dy^2=x(x+3)(x-1).
$$
It should be noted that $E$ and $E^d$ are $\Q(\sqrt{d})$-isomorphic.

\begin{corollary}
Let $d$ be a squarefree integer. Then there is a non-constant arithmetic progression of four squares over $\Q(\sqrt{d})$ if and only if $\rank{E^d}\ne 0$; in this case, there exist infinitely many such progressions.
\end{corollary}

\begin{proof} We are going to compute the structure of $E(\Q(\sqrt{d}))$. Since the $2$-torsion subgroup of $E$ is defined over $\Q$, by applying Kwon's results \cite{kwon} we see that $E(\Q(\sqrt{d}))_{\tors}$ and $E(\Q)_{\tors}$ are equal. Thus Proposition \ref{prop1} shows that if there exists a non-constant arithmetic progression of four squares over $\Q(\sqrt{d})$, then there exist infinitely many or, equivalently, $\rank{E(\Q(\sqrt{d}))}\ne 0$. Now, since 
\begin{equation}\label{rank-formula}
\rank E(\Q(\sqrt{d})) = \rank E(\Q)+ \rank E^d(\Q),
\end{equation}
the statement follows from $\rank E(\Q)=0$.
\end{proof}

\vspace{.3cm}

Therefore, the problem to decide for which quadratic fields $\Q(\sqrt{d})$ there exist non-constant arithmetic progressions of four squares is reduced to the question whether the rank of $E^d(\Q)$ is positive or not. In Section \ref{section-explicit} we give a partial solution to this problem and for some quadratic number fields we also give explicit four term arithmetic progressions of squares.

\vspace{.2cm}

For number fields of higher degree we have the following result:

\begin{theorem}
There are infinitely many cyclic cubic number fields in each of which there exist infinitely many non-constant arithmetic progressions of four squares.
\end{theorem}

\begin{proof}
Applying \cite[Theorem 6.1]{FKK} to the elliptic curve $E=X_0(24)$, we find $\#E(\Q)=8>6$. Hence, for infinitely many cyclic cubic extensions $K/\Q$ we have $\rank E(K) > \rank E(\Q)=0$.
\end{proof}

\section{$\theta$-congruent numbers}\label{section-congruent}

A positive integer $n$ is called a congruent number if there exists a right triangle with rational sides and area equal to $n$. The problem to decide whether a given integer is a congruent number has been studied since Diophantus. In 1983, Tunnell \cite{tunnell} found a deterministic criterion for this problem; if the Birch \& Swinnerton--Dyer conjecture is true, then his criterion also can be used to determine congruent numbers. The notion of congruent numbers was extended by Fujiwara \cite{fujiwara} to rational $\theta$-triangles, that are triangles with rational sides where one angle is equal to $\theta$. Note that for such a triangle,  $\cos \theta=s/r$ for some coprime integers $r,s$ with $r>0$. It follows that $\sin\theta=\alpha_{\theta}/r$, where $\alpha_{\theta}:=\sqrt{r^2-s^2}$ is uniquely determined by $\theta$. Then $\theta$-congruent numbers are defined as follows:

\begin{definition}
Let be $\theta \in [0,\pi)$. A positive integer $n$ is a $\theta$-congruent number if there exists a rational $\theta$-triangle whose area is equal to $n\alpha_\theta$.
\end{definition}

Therefore, $\pi/2$-congruent numbers coincide with the ordinary congruent numbers (in which case $r=1$ and $s=0$). Generalizing the case of ordinary congruent numbers, there is a characterization of $\theta$-congruent numbers in terms of rational points on elliptic curves.

\begin{theorem}{\cite{fujiwara,fujiwara2}}
For $\theta \in [0,\pi)$ and $n\in\mathbb{N}$, define the elliptic curve
$$
E_{n,\theta} : y^2=x(x+(r+s)n)(x-(r-s)n).
$$
Then, $n$ is a $\theta$-congruent number if and only if there exists a rational point on $E_{n,\theta}$ of order greater than $2$. Moreover, if $n\ne 1,2,3,6$, then $n$ is a $\theta$-congruent number if and only if $\rank E_{n,\theta}(\Q)>0$.
\end{theorem}

\noindent {\bf Remark:} Note that the elliptic curve $E_{n,\theta}$ is the $n$-twist of $E_{1,\theta}$. Therefore, whenever $n\ne 1,2,3,6$, to prove that $n$ is a $\theta$-congruent number is equivalent to show that the rank of the $n$-twist of the elliptic curve $E_{1,\theta}$ is non-zero. Another interesting remark is that $E_{n,\pi-\theta}$ is the $(-n)$-twist of $E_{1,\theta}$ .
\vspace{.3cm}

Several papers \cite{fujiwara,fujiwara2,goto,HK,kan, yoshida,yoshida2} have been studying the $\theta$-congruent number problem for $\theta\ne \pi/2$. For our purposes the cases $\theta=\pi/3$ and $\theta=2\pi/3$ are of special interest (see Section \ref{section-explicit}). In these cases, the curve $E_{1,\pi/3}$ is the curve \verb+24A1+ and $E_{1,2\pi/3}$ is the curve \verb+48A1+ in Cremona's tables, respectively.

The following table resumes all known results on $\pi/3$-congruent and $2\pi/3$-congruent numbers (see \cite{fujiwara,HK,kan, yoshida,yoshida2}). Here '?' indicates that in this case we do not know whether $n$ is $\theta$-congruent or not.

{\footnotesize
\begin{center}
\begin{table}[!h]
\begin{tabular}{|c||c|c|c|c||c|c|c|c|}
\hline
$p\ge 5$ prime & \multicolumn{4}{|c||}{Is $n$ a $\pi/3$-congruent number?} & \multicolumn{4}{|c|}{Is $n$ a $2\pi/3$-congruent number?}\\
\hline
$p\bmod\,24$ & $n=p$ & $n=2p$ &$n=3p$ &$n=6p$ & $n=p$ & $n=2p$ &$n=3p$ &$n=6p$\\
\hline
 $1$ & $?$ & $?$ & $?$ & $?$ & $?$ &  $?$ & $?$ & $?$ \\
\hline
 $5$ & no & $?$ & no & $?$ & $?$ &  $?$ & $?$ & no \\
\hline
 $7$ & no & no &  $?$ & $?$ & no &  $?$ & $?$ & no \\
\hline
 $11$ & $?$ & $?$ & no & $?$ & no &  $?$ & $?$ & $?$ \\
\hline
 $13$ & $?$ & no &  $?$ & $?$ & no &  no & $?$ & no \\
\hline
 $17$ & $?$ & $?$ &  no & $?$ & $?$ &  $?$ & no & no \\
\hline
 $19$ & no & $?$ &  no & $?$ & $?$ & no & $?$ & $?$\\
\hline
 $23$ & yes & yes &  $?$ & yes & yes &  yes & yes & $?$ \\
\hline
\end{tabular}\\[2mm]
\caption{}\label{table2}
\end{table}
\end{center}
}
\vspace{-8mm}

Let $n$ be a squarefree positive integer from one of the residue classes $1,7$ or $13\bmod\,24$. Yoshida \cite{yoshida,yoshida2} proved along the lines of Tunnell's approach to the congruent number problem that $n$ is not a $2\pi/3$-congruent number if the number of representations of $n$ by the ternary quadratic forms 
$$
X^2+3Y^2+144Z^2\qquad\mbox{and}\qquad 3X^2+9Y^2+16Z^2
$$
with integral $X,Y,Z$ are not equal; if the conjecture of Birch and Swinnerton-Dyer is true, then also the converse implication holds, i.e., $n$ is $2\pi/3$-congruent if the number of representations are identical. In particular, it follows from the theory of quadratic forms that all primes $p\equiv 7$ or $13\bmod\,24$ are not $2\pi/3$-congruent numbers. Moreover, if $n\neq 1$ is squarefree with $n\equiv 1,7$ or $19\bmod\,24$, then he showed that $n$ is not a $\pi/3$-congruent number if the number of integer representations of $n$ by the ternary quadratic forms
$$
X^2+12Y^2+15Z^2+12YZ\qquad\mbox{and}\qquad 3X^2+4Y^2+13Z^2+4YZ
$$
are not equal; if the conjecture of Birch and Swinnerton-Dyer is true, then also the converse implication is true. Here it follows that no prime $p\equiv 7\bmod\,24$ is a $\pi/3$-congruent number. For other residue classes Yoshida obtained analogous statements with, of course, different quadratic forms, which explain the 'no's in the above table. Note that all affirmative 'yes' rely on primes $p\equiv 23\bmod\,24$ by a theorem of Kan \cite{kan}. We collect these results on primes in the following theorem; for the other cases we refer to the mentioned papers.

\begin{theorem} \cite{kan,yoshida,yoshida2}
There is no prime number $p\equiv 7,11,13\bmod\,24$ which is a $2\pi/3$-congruent number, and there is no prime number $p\equiv 5,7,19\bmod\,24$ which is a $\pi/3$-congruent number where $p>5$. On the contrary, any prime $p\equiv 23\bmod\,24$ is a $\theta$-congruent number for both $\theta=\pi/3$ and $\theta=2\pi/3$.
\end{theorem}

A classical conjecture, supported by numerical evidence (based on computations from the 1970s), states that all squarefree positive integers congruent $5,6$ or $7$ modulo $8$ are congruent numbers. In fact, this conjecture is a direct consequence of the Birch \& Swinnerton--Dyer conjecture. Similar to this conjecture the following one covers the cases of $\pi/3$- and $2 \pi/3$-congruent numbers, respectively.
\vspace{.3cm}

\begin{conjecture}\label{conjecture}
Let $n$ be a squarefree positive integer.
\begin{itemize}
\item If $n\equiv\,11,13,17,23\,\mbox{mod $24$}$, then $n$ is a $\pi/3$-congruent number.
\item If $n\equiv\,5,17,19,23\,\mbox{mod $24$}$, then $n$ is a $2\pi/3$-congruent number.
\end{itemize}
\end{conjecture}

\section{Euler's concordant forms}\label{section-concordant}

Another equivalent formulation of the congruent number problem is the following: a positive integer $n$ is a congruent number if and only if the system of diophantine equations
$$
\left\{
\begin{array}{c}
x^2 + ny^2=t^2\\
x^2 - ny^2=z^2
\end{array}
\right.
$$
has a solution $x,y,z,t\in\Z$ with $xy\ne 0$. In 1780, Euler \cite{euler} gave another generalization of the congruent number problem. He was interested to classify those pairs of distinct non-zero integers $M$ and $N$ for which there exist $x,y,z,t\in\Z$ with $xy\ne 0$ such that
$$
\left\{
\begin{array}{c}
x^2+My^2=t^2\\
x^2+Ny^2=z^2
\end{array}
\right.
$$
This is known as Euler's concordant forms problem. If the above diophantine system has a solution, then the pair $(M,N)$ is said to be concordant, otherwise discordant. In the particular case when $M=-N$ this yields the congruent number problem. As for the congruent number problem and its generalization to $\theta$-congruent numbers, there is a characterization due to Ono \cite{ono} for Euler's concordant forms problem in terms of rational points of an elliptic curve. 

\begin{theorem}{\cite{ono}}
For $M,N\in\Z$ such that $NM(M-N)\ne 0$, define the elliptic curve
$$
E_{M,N}:y^2=x(x+M)(x+N).
$$
Then, the pair $(M,N)$ is concordant if and only if there exists a rational point on $E_{M,N}$ of order $\neq 1,2$ or $4$. In particular, if $\rank E_{M,N}(\Q)$ is positive, then $(M,N)$ is concordant.
\end{theorem}

Using Waldspurger's results and Shimura's correspondence {\it a la} Tunnell, Ono obtained several results on the ranks of twists of $E_{M,N}$. In particular, if the number of integer representations of an odd positive squarefree integer $n$ by the ternary quadratic forms
$$
X^2+2Y^2+12Z^2 \qquad\mbox{and}\qquad 2X^2+3Y^2+4Z^2
$$
is not equal, then $(M,N)=(6n,-18n)$ is discordant; if the conjecture of Birch and Swinnerton-Dyer is true, then also the converse implication is true. Moreover, for $r$ being an odd integer with $1\leq r\leq 24$, he showed that there are infinitely many positive squarefree integers $n\equiv\,r\bmod\,24$ such that
\begin{itemize}
\item $(M,N)=(6n,-18n)$ is discordant.
\item $(M,N)=(9n,-3n)$ is discordant, where $r\ne 7,15,23$.
\end{itemize}

\section{Four squares in arithmetic progressions over quadratic fields}\label{section-explicit}

The existence of a non-constant four term arithmetic progression of squares over a quadratic field is determined by the rank of twist of the elliptic curve $E=X_0(24)$. For a real-quadratic field $\Q(\sqrt{d})$ the elliptic curve $E^d$ is equal to $E_{d,\pi/3}$, whereas for an imaginary-quadratic field $\Q(\sqrt{-d})$ it is $E_{d,2\pi/3}$. 

Thus, we may use the information of the table \ref{table2} to deduce information about the rank of $E^d(\Q)$ from $E_{d,\pi/3}$ and $E_{d,2\pi/3}$, respectively, corresponding to $d$ being positive or negative. This proves the information provided by table \ref{table1} from the introduction. Moreover assuming the conjecture \ref{conjecture} we have that there exists a non-constant arithmetic progression of four squares over $\mathbb{Q}(\sqrt{d})$ if $d\equiv\,11,13,17,23\,\mbox{mod $24$}$ in the real case and $d\equiv\,1,5,7,19\,\mbox{mod $24$}$ in the imaginary case.

Further, note that for $(M,N)=(-1,3)$ or $(3,-1)$ we have $E_{d,\pi/3}=E_{M,N}$; moreover, if $(M,N)=(1,-3)$ or $(-3,1)$, then $E_{d,2\pi/3}=E_{M,N}$. In these cases we can use Ono's results \cite{ono} on the rank of $E_{M,N}(\Q)$.

\subsection{Explicit examples}

Using \verb+MAGMA+, we may compute the rank of $E^d(\Q)$; if this rank is positive, using the parametrization from Theorem \ref{parametrization}, we can also compute an explicit arithmetic progression of four squares. The following two tables list explicit examples of such progressions according to $\Q(\sqrt{d})$ being an imaginary-quadratic or a real-quadratic number field for the range $|d|\le 40$. Each of the tables consists of three colums; the first column indicates the value of $d$, the second and third one give an example for $a,r\in\Q(\sqrt{d})$ such that $a^2,a^2+r,a^2+2r,a^2+3r$ forms an arithmetic progressions of squares over $\Q(\sqrt{d})$, where $\alpha:=\sqrt{d}$.
\medskip

\begin{center}
\begin{tabular}{|c|c|c|}
\hline
$d$ & $a$ & $r$\\
\hline
$-5$ & $(-4 \alpha - 73)/2$ & $42 \alpha - 840$  \\ 
$-10$ & $(-60 \alpha - 629)/2$ & $-14322 \alpha - 55440$  \\ 
$-14$ & $-6 \alpha - 361$ & $8580 \alpha - 65520$  \\ 
$-15$ & $(\alpha - 19)/4$ & $\alpha - 15$  \\ 
$-17$ & $3339440 \alpha - 57973177$ & $219447603254880 \alpha - 2180390987174400$  \\ 
$-19$ & $(-26589360 \alpha + 21015523)/2$ & $22024049983320 \alpha + 2196495218332800$  \\ 
$-21$ & $-160 \alpha - 2393$ & $-590920 \alpha - 3141600$  \\ 
$-22$ & $-1224720 \alpha - 2179673$ & $-3235062111120 \alpha + 19986461510400$  \\ 
$-23$ & $(2625 \alpha + 11951)/4$ & $-1231230 \alpha + 6592950$  \\ 
$-29$ & $(22143940 \alpha - 361130617)/2$ & $3478556113902870 \alpha - 13110939890248200$  \\ 
$-33$ & $-612000 \alpha + 1945781$ & $673901512480 \alpha + 8195655283200$  \\ 
$-34$ & $319440 \alpha + 650807$ & $-95395404720 \alpha + 2288266041600$  \\ 
$-39$ & $(36 \alpha - 683)/2$ & $10450 \alpha - 51480$  \\ 
\hline
\end{tabular}
\end{center}
\vspace{.1cm}
\begin{center}
\begin{tabular}{|c|c|c|}
\hline
$d$ & $a$ & $r$\\
\hline
$6$ & $(2 \alpha + 1)/2$ & $25 \alpha + 60$  \\ 
$10$ & $(10 \alpha - 19)/2$ & $33 \alpha - 60$  \\ 
$11$ & $(-1320 \alpha + 2843)/2$ & $-7242060 \alpha + 23839200$  \\ 
$13$ & $1440 \alpha + 5183$ & $1323960 \alpha + 4773600$  \\ 
$17$ & $-15 \alpha - 1511$ & $-555360 \alpha + 1591200$  \\ 
$21$ & $-36 \alpha + 163$ & $2200 \alpha - 10080$  \\ 
$22$ & $(-750 \alpha + 3529)/2$ & $453705 \alpha - 2009700$  \\ 
$23$ & $19476668640 \alpha + 90283636367$ & 
\begin{tabular}{c}
$\!\!\!\!\!\!\!\!\!\!\!\!\!1725783576049531078080 \alpha +$\\ 
$\qquad\qquad 8274631385821773772800$
\end{tabular}\\ 
$30$ & $(18 \alpha + 19)/2$ & $-413 \alpha + 1260$  \\ 
$34$ & $6860 \alpha - 12239$ & $-12575640 \alpha + 43982400$  \\ 
$35$ & $-84 \alpha + 487$ & $-28968 \alpha + 171360$  \\ 
$37$ & 
$38306628360 \alpha - 276487794001$ & 
\begin{tabular}{c}
$\!\!\!\!\!\!\!\!\!\!\!\!\!10021678513795431723240 \alpha -$\\
 $\qquad\qquad 24114970612028472976800$  
\end{tabular}\\ 
$39$ & $720 \alpha + 3869$ & $7990640 \alpha + 49795200$  \\ 
\hline
\end{tabular}
\end{center}

\vspace{.3cm}

\subsection{Pythagorean triples}

Next we use Pythagorean triples to construct arithmetic progressions of four squares over quadratic number fields. In the following cases the arithmetic progression consists of four integers, all of them being a square over a specific quadratic number field.

It is well-known that for arbitrary $a,b\in \Z$ the triple $(a^2-b^2,2ab,a^2+b^2)$ defines a Pythagorean triple. Then $(a^2+b^2)^2-4n, (a^2+b^2)^2,(a^2+b^2)^2+4n$
forms an arithmetic progression of three squares where $n=ab(a^2-b^2)$. If we add a new square term, $\alpha^2=(a^2+b^2)^2+8ab(a^2-b^2)$ say, we obtain an arithmetic progression of four squares over $\Q(\sqrt{(a^2+b^2)^2+8ab(a^2-b^2)})$. In order to construct a quadratic number field $\Q(\sqrt{d})$, where $d$ is a squarefree integer, we define the Thue equation
$$
F(x,y)=(x^2+y^2)^2+8xy(x^2-y^2)=d.
$$
Of course, here we are interested in the set of integer solutions. Using $\verb+MAGMA+$, we have observed that for $|d|<100$ there exist integer solutions in the cases $d=-71, -47, -23, 73$, all of them being congruent to $1\bmod\,24$. By the above construction this yields the following arithmetic progressions:

\vspace{.1cm}
\begin{center}
\begin{tabular}{|c|c|c|}
\hline
$d$ & Four square in arithmetic progression over $\Q(\sqrt{d})$\\[.1cm]
\hline
$-71$ & $(\sqrt{-71})^2,7^2,13^2,17^2$\\[.1cm]
$-47$ & $(\sqrt{-47})^2,17^2,25^2,31^2$\\[.1cm]
$-23$ & $(\sqrt{-23})^2,1,5^2,7^2$\\[.1cm]
$73$ & $1,5^2,7^2,(\sqrt{73})^2$\\[.1cm]
\hline
\end{tabular}
\end{center}

It is straightforward to compute that $(a^2+b^2)^2+8ab(a^2-b^2) \equiv\,1\bmod\,24$ for $a,b\in \Z$ with coprime $a,b$ and $a\not\equiv\,b\bmod\,2$. Therefore, our construction is restricted to the case of number fields $\Q(\sqrt{d})$ with $d\equiv\,1\bmod\,24$. Moreover, all arithmetic progressions discovered by this method satisfy that the difference of any two successive members is divisible by $24$. 

\section{Average results}

We conclude by discussing some average results for the central values of $L$-functions associated with elliptic curves. Let $E$ be an elliptic curve over $\Q$ and denote by $L(E,s)$ its elliptic curve $L$-function. Roughly speaking, the yet unproved Birch \& Swinnerton--Dyer conjecture states that the order of non-vanishing of $L(E,s)$ is equal to the rank of $E$. Kolyvagin \cite{koly} has shown that if $E$ is a modular elliptic curve with $L(E,1)\neq 0$, then the rank of $E$ is equal to zero. By the proof of Wiles et al. \cite{wiles,wiles2} of the Shimura-Taniyama conjecture any elliptic curve over $\Q$ is modular, however, the corresponding statement for quadratic number fields is not known to be true. Note that
$$
L(E(\Q(\sqrt{d})),s)=L(E,s)L(E^d,s),
$$
where $E^d$ denotes the quadratic twist of $E$ by $d$; here it suffices to consider squarefree integers $d$. This formula directly corresponds to (\ref{rank-formula}). Hence in our case, in order to have rank zero for $E(\Q(\sqrt{d}))$ we need the non-vanishing of $L(E^d,1)$ since $L(E,1)=0.53...$; this computations has been made using \verb+SAGE+.

We may ask for the statistical behaviour as $d$ varies. Goldfeld \cite{gold} has conjectured that a positive proportion of $0<\vert d\vert \leq X$ have the property that $L(E^d,1)$ is non-vanishing. This has been established only in exceptional cases. For instance, Heath--Brown \cite{hb} confirmed this conjecture for the congruent number elliptic curve. Moreover, Ono \& Skinner \cite[Corollary 2]{onoskin} have proved that if $E$ is an elliptic curve over $\Q$ with conductor $\leq 100$, then either $E^{-p}$ or $E^p$ has rank zero for a positive proportion of primes $p$. In the special case of our elliptic curve we obtain:

\begin{corollary}
For a positive proportion of primes $p$, there are either no non-constant arithmetic progressions of four squares in $\Q(\sqrt{-p})$ or $\Q(\sqrt{p})$.
\end{corollary}

If the Birch \& Swinnerton--Dyer conjecture is true, then the vanishing of $L(E^d,1)$ would imply that the rank of $E(\Q(\sqrt{d})$ is positive, and so from (\ref{rank-formula}) and Proposition \ref{prop1} it would follow that there exist infinitely many non-constant arithmetic progressions of four squares over $\Q(\sqrt{d})$. M.R. Murty \& V.K. Murty \cite{mm} have shown that for $L(E,1)\neq 0$ there are infinitely many fundamental discriminants $d<0$ such that $L(E^d,s)$ has a simple zero at $s=1$; this result was independently obtained by Bump, Friedberg \& Hoffstein \cite{bump}. Hence, in our special case we may deduce that there exist infinitely many imaginary quadratic fields $\Q(\sqrt{d})$ each of which containing infinitely many non-constant arithmetic progressions of four squares subject to the truth of the Birch \& Swinnerton--Dyer conjecture.

%%%%%%%%%%%%%%%%%%%%%%%%%%%%%%%%%%%%%%%%%%%%%%%%%%
%%%%%%%%%%%%%%%%%%%%%%%%%%%%%%%%%%%%%%%%%%%%%%%%%%
%%%%%%%%%%%%%%%%%%%%%%%%%%%%%%%%%%%%%%%%%%%%%%%%%%

\end{document}